\documentclass[twoside]{amsart}
\usepackage{pb-diagram,epic,eepic}

\theoremstyle{plain}
\newtheorem{thm}{Theorem}[section]

\author{Valery Alexeev} 
\address{Department of Mathematics\\
University of Georgia\\
Athens, GA 30602}
\email{valery@math.uga.edu}

\newcommand{\bL}{{\mathbb L}}

\newcommand{\bP}{{\mathbb P}}

\newcommand{\bZ}{{\mathbb Z}}

\newcommand{\bR}{{\mathbb R}}
\newcommand{\bC}{{\mathbb C}}

\newcommand{\cO}{{\mathcal O}}

\newcommand{\inv}{^{-1}}

\newcommand{\hbL}{\widehat \bL}

\newcommand{\wP}{\widetilde P}

\newcommand{\wL}{\widetilde L}

\newcommand{\wDelta}{\widetilde\Delta}

\DeclareMathSymbol{\curvearrowright}{\mathrel}{AMSb}{"79}
\DeclareMathSymbol\rightsquigarrow {\mathrel}{AMSa}{"20}
\DeclareMathSymbol\square {\mathord}{AMSa}{"03}
\DeclareMathSymbol{\ltimes}         {\mathbin}{AMSb}{"6E}
\DeclareMathSymbol{\nmid}           {\mathrel}{AMSb}{"2D}
\DeclareMathSymbol{\twoheadrightarrow}  {\mathrel}{AMSa}{"10}

\newcommand{\ratmap}{- \kern -3pt \to}

\newcommand{\acts}{\curvearrowright}

\newcommand{\GL}{\operatorname{GL}}

\newcommand{\Pic}{\operatorname{Pic}}

\newcommand{\Sym}{\operatorname{Sym}}

\theoremstyle{plain}

\newtheorem{conj}[thm]{Conjecture}

\theoremstyle{definition}

\newtheorem{defn}[thm]{Definition}

\newtheorem{exmp}[thm]{Example}

\theoremstyle{remark}

\newcommand{\XR}{X_{\bR}}
\newcommand{\Atau}{\overline{\operatorname{A}}_g^{\tau}}
\newcommand{\Avor}{\overline{\operatorname{A}}_g^{\text{Vor}}}
\newcommand{\Jac}{\operatorname{Jac}}
\newcommand{\mom}{\operatorname{Mom}}
\newcommand{\str}{\operatorname{Stratum}}
\newcommand{\vor}{\operatorname{Vor}}
\newcommand{\Del}{\operatorname{Del}}
\newcommand{\APg}{\overline{\operatorname{AP}}_g}
\newcommand{\bld}{\widehat\bL_{\delta,\bR}}
\newcommand{\bl}{\widehat\bL}

\title{On extra components in the functorial compactification of
  $A_g$} 
\date{April 25, 1999}
\begin{document}
\bibliographystyle{amsalpha+}
\maketitle

\section*{}
Recall the following from the theory of toroidal compactifications of
moduli of polarized abelian varieties (Mumford et al \cite{AMRT} over
$\bC$, Faltings and Chai \cite{FaltingsChai90} over $\bZ$). Denote
$X=\bZ^g$ and let $C$ be the convex hull in the space
$\Sym^2(X_{\bR}^*)$ of semipositive symmetric matrices $q$ with
rational null-space. For any \emph{admissible} $\GL(X)$-invariant
decomposition $\tau$ of $C$ (i.e.  it is a face-fitting decomposition
into finitely generated rational cones such that there are only
finitely many cones modulo $\GL(X)$) there is a compactification
$\Atau$ of the moduli space $A_g$ of principally polarized abelian
varieties. $\Atau$ comes with a natural stratification, and strata
correspond in a 1-to-1 way to cones in $\tau$ modulo $\GL(X)$.  There
are infinitely many such decompositions $\tau$ and none of the seems
to better that another.  True, some decompositions are smooth and
projective but still there are infinitely many of these as well.

There is, however, a decomposition $\tau_{\vor}$ for the 2nd Voronoi
decomposition which has a nicer geometric description. Strata of this
compactification $\Avor$ still are in a bijection with cones of
$\tau_{\vor}$, however, these cones now correspond in a 1-to-1 way to
\emph{special} $X$-periodic face fitting decompositions of $X_{\bR}$
with vertices in $X$, called \emph{Delaunay}.  A form $q$ defines a
distance function $d_q$ on $X_{\bR}$ and a cell of the Delaunay
decomposition $\Del_q$ (Delaunay cell) is a convex hull of integral
points circumscribed by an ``empty sphere'', a sphere that does not
contain any integral points in its interior. As an example, for $g=2$
there are 4 such decompositions: by 4-gons, by triangles, by infinite
strips and finally the decomposition consisting of one big cell
covering the whole plane.  The decompositions appearing are either
polytopal or the preimages of such from a lower dimension.

On the other hand, in \cite{Alexeev_CMAV} I have constructed a
functorial compactification $\APg$ of $A_g$ as the moduli of triples
$G\acts P\supset \Theta$ whose geometric fibers have the following
description: $G$ is semiabelian, $P$ is projective reduced connected
and $\Theta$ is a Cartier divisor, all satisfying a few natural
conditions:
\begin{enumerate}
\item $P$ is seminormal,
\item there are only finitely many orbits,
\item $\Theta$ does not contain any orbit entirely,
\item for any $p\in P$, the stabilizer of $p$ is connected and
  reduced and lies in the toric part of $G$. 
\end{enumerate}

This compactification comes with a
stratification as well, with strata corresponding to \emph{all}
$X$-periodic face fitting decompositions of $X_{\bR}$ with vertices in
$X$. Just from this rough description we see that $\Avor$ and $\APg$
must be very closely related and at the same time be different. The
connection, according to \cite{Alexeev_CMAV}, is that $\Avor$
coincides with the main irreducible component of $\APg$. The most
obvious difference is that, unlike $\Avor$, $\APg$ $\,$ may have Extra
Types of irreducible components, ETs for short, which will be the
focus of this note. We would like to discuss:
\begin{enumerate}
\item Where and why ETs appear.
\item How to find ETs and how to study them.
\item Some concrete evidence of ETs.
\item Dimension 4 case in detail.
\item Relationship between ETs and the Jacobian locus. 
\end{enumerate}

\section{Where and why ETs appear}

There is actually a very simple reason for their appearance. Say,
$(P',\Theta')$ is a pair with an abelian action and let $(P,\Theta)$
be its degeneration. Assume that $(P,\Theta)$ has several components
$P_i$ and set $\Theta_i=\Theta|_{P_i}$. Think of the deformations of
$(P,\Theta)$ that keep the decomposition into the irreducible
components. These deformations correspond to deformations of components
$(P_i,\Theta_i)$ that are compatible on intersections. In the simplest
case, when the intersections are elementary, there are, perhaps, no
gluing conditions at all. Now, it is not hard to imagine the situation
when the sum of the dimensions of the deformation spaces for pairs
$(P_i,\Theta_i)$ is larger than the dimension of the deformation space
of the constrained smooth pair $(P',\Theta')$, i.e. $g(g+1)/2$. In
this case there must be another irreducible component.

So, the answer to the ``where'' part is that ETs appear at the
boundary of the world as we know it.

\section{How to find ETs and how to study them}

There are two basic methods:
\begin{enumerate}
\item Find a periodic non-Delaunay decomposition. The first such
  decomposition appears in dimension 4, and we will describe it in
  more detail below.
\item For a Delaunay decomposition $\wDelta$, compute the dimension
  of $\str\wDelta$ in $\APg$. If it is higher than dimension of the
  corresponding stratum in $\Avor$ then there must be an ET nearby.
  There is a remarkably simple formula for the first of these
  dimensions which we are now going to describe.
\end{enumerate}

Let $\wDelta$ be an $X$-periodic face-fitting decomposition of $\XR$
with vertices in $X$. We will identify $\wDelta$ with its quotient
$\Delta$ which is a decomposition of the real torus $\XR/X\simeq
\bR^g/\bZ^g$.  If we work over $\bC$, this decomposition is directly
related to the properties of each pair $(P,\Theta)$ in the
$\Delta$-stratum of $\APg$, and this connection is described by using
the \emph{moment map}.  (Even if the base field is not $\bC$, this
decomposition describes the main properties of $(P,\Theta)$ very
faithfully.)

Say, $G\acts P\supset \Theta$ is a triple as before, and start with
the case when $G=T=(\bC^*)^g$ is a torus with the character group $X$.
The moment map sends $P(\bC)$ to its quotient by the action of the
compact torus $CT=(S^1)^g = U(1)^g \subset T$. If $L=\cO(\Theta)$ is
$T$-linearized then we can describe the moment map more directly. Let
$\theta\in H^0(P,L)$ be an equation of $\Theta$. $H^0(P,L)$ splits
into the direct sum of $T$-eigenspaces and we can write $\theta$ as a
finite sum $\sum_{x\in X} \xi_x$. Then
\begin{displaymath}
  \mom: P(\bC)\ni p\mapsto \frac{\sum |\xi_x(p)|^2\cdot x}{\sum
    |\xi_x(p)|^2} \in X_{\bR}
\end{displaymath}
It is well defined if at each $p\in P$ at least one of $\xi_x$'s is
not zero, and that is one of the conditions on $(G,P,\Theta)$ that is
satisfied by our definition of a triple.

If $L$ is not linearized then $(P,L)$ is in a canonical way the
quotient by $X$ of a pair $(\wP,\wL)$ with linearized $\wL$ (the
scheme $\wP$ is only locally of finite type but the action of $X$ is
properly discontinuous in Zariski topology). We take $\tilde\theta\in
H^0(\wP,\wL)$ to be the pullback of $\theta$.  It is a fact that for
any $\tilde p\in\wP$ the sum $\sum_{x\in X} \xi_x$ is finite as almost
all of $\xi_x$ vanish at $\tilde p$. The moment map
$\mom_{\wP,\tilde\theta}$ commutes with the translation action of $X$,
and we define $\mom_{P,\theta}$ from the following diagram
\begin{displaymath}
    \begin{diagram}
    \node{\wP(\bC)}\arrow{e,t}{\mom_{\wP,\tilde\theta}}
    \arrow{s,l}{/X} 
    \node{\XR}\arrow{s,r}{/X} \\
    \node{P(\bC)}\arrow{e,t}{\mom_{P,\theta}}
    \node{\XR/X}
  \end{diagram}
\end{displaymath}

Finally, if $G$ is an arbitrary semiabelian variety then we can (there
is a choice involved) write $G$ as the quotient $(\bC^*)^g/X'$,
$X\supset X'=\bZ^{g'}$ and $g'=$ the dimension of the abelian part of
$G$, making $P$ into a ``toric'' variety. The $T$-action here is not
algebraic, of course, but we can still repeat the above definition.
The sum $\sum |\xi_x(p)|^2\cdot x$ in this case is truly infinite,
however, it is convergent because theta functions have exponential
decline for large $x$.

With these definitions in mind, let $(G,P,\Theta)$ be a triple
corresponding to a point in $\str[\Delta]\subset\APg$. Then the moment
map sends $P$ to $\XR/X$ and for any $z\in \delta^0$, interior of
$\delta$, one has $\mom\inv(z)=(S^1)^{\dim\delta}$. Moreover, $P$ has
as many irreducible components as there are cells in $\Delta$, and
these intersect exactly in the way the corresponding cells intersect.
Moreover, if $\delta$ is a polytope then the irreducible component
$P_{\delta}$ is the projective toric variety corresponding to
$\delta$, glued in a way that the decomposition suggests (in
particular, all of its $0$-dimensional orbits are glued together).

As triples degenerate, so do the moment maps. If the fiber over a
particular point was $0$-dimensional, it is not going to get bigger.
On the other hand, a big fiber may get ``squashed'' to a smaller one
in the limit, so the degeneration of triples of type $\Delta$ must
correspond to a subdivision of $\Delta$. A familiar example is that of
an elliptic curve degenerating to a nodal curve. An elliptic curve has
a moment map to a circle whose every fiber is $S^1$, and in the limit
a fiber over one point becomes $0$-dimensional. An elliptic curve
corresponds to a decomposition $\wDelta$ with one big cell, and the
nodal curve -- to a decomposition of $\bR$ into intervals. In higher
dimensions exactly the same happens but the decompositions $\Delta$
that appear are more sophisticated.

\begin{defn}
  For a ``simple'' cell $\delta$ we define $\bld$ to be the space of
  all $\bR$-valued functions on $\delta\cap X$ modulo the subspace of
  linear nonhomogeneous functions. $C^0(\Delta,\bl)$ is defined as
  \begin{displaymath}
    \oplus_{\dim\delta=g} \, \bld
  \end{displaymath}
  and $H^0(\Delta,\bl)$ as the subspace of $C^0$ consisting of
  functions that coincide on ``intersections''. 
\end{defn}

The reader will notice that we are computing here the space of global
sections of a certain constructible sheaf $\bl$ on $\XR/X$ which is
constant on locally closed strata defined by $\Delta$.

\begin{thm}[\cite{Alexeev_CMAV}]
  $\dim\str[\Delta] = h^0(\Delta,\bl)$
\end{thm}

The meaning of ``simple'' and ``intersections'' should be clear from
the following instructive examples.

\begin{exmp}
  Say, $\wDelta$ consists of just one big cell, $\bR^g$ -- this is the
  stratum corresponding to triples with abelian $G$. Then $C^0(\bl)$
  is the space of all real-valued functions on $X=\bZ^g$ modulo the
  $(g+1)$ dimensional subspace of linear nonhomogeneous functions.
  $H^0(\bl)$ is the subspace of functions that are invariant under the
  translations by $X$. If $[f]$ is the equivalence class function of
  such a function then for any $y\in X$ the function
  $g(x)=f(y+x)-f(x)$ must be linear, so $f(x)$ must be quadratic.
  Therefore, $H^0(\bl)$ is the space of quadratic modulo linear
  functions and its dimension is $g(g+1)/2$, i.e. the dimension of
  $A_g$, as expected.
\end{exmp}

\begin{exmp}
  Let $\wDelta$ be the decomposition of $\bR^2$ into squares. There is 
  only 1 maximal-dimensional cell modulo $X$, a square. It is easy to
  see that for any polytopal cell appearing in an $X$-periodic
  decomposition one has
  \begin{displaymath}
    l_{\delta}= \dim\bld=
    \#(\text{vertices of } \delta) - \dim\delta -1
  \end{displaymath}
  For the square we have $l=1$ and for the intervals on the boundary:
  $l=0$, so there are no gluing conditions. Hence, the dimension of
  the corresponding stratum is $h^0(\bL)=1$.
\end{exmp}

\begin{exmp}\label{exmp:decomps_with_simplicial_cells}
  Let us generalize the previous example slightly. Assume that all
  cells of $\wDelta$ are polytopal and that all maximal cells are
  simplicial, i.e. all of its proper faces are simplices.  Then again
  $C^1(\bl)=\oplus_{\dim\delta=g-1} \, \bld=0$ and one obtains
  \begin{eqnarray}\label{eqn:main}
    h^0(\hbL) = \sum_{\dim\delta=g} l_{\delta}. 
  \end{eqnarray}
This very simple formula has many applications as we will see.
\end{exmp}

\cite{Alexeev_CMAV} contains a more detailed information about
$\str[\Delta]$. In particular, if all the cells are polytopal then the
normalization of the closure of this stratum is the quotient by the
finite group $\Sym\Delta$ of the projective toric variety
corresponding to the so called generalized secondary polytope
$\Sigma(\Delta)$.

\section{Some concrete evidence of ETs}

\begin{thm}\label{thm:dim_le_three}
  There are no ETs in dimension $g\le3$.
\end{thm}

\begin{thm}\label{thm:dim_four}
  For $g=4$ there is exactly 1 ET and it is isomorphic to $\bP^2$.
\end{thm}

\begin{thm}\label{thm:dim_large}
  For $g\ge5$, $2^g-g-3\le \dim\APg <g!$
\end{thm}
Hence, the maximal dimension of ETs grows at least exponentially. The
dimension of the ``main'' component, of course, grows as a particular
polynomial of degree 2, namely $g(g+1)/2$. One might say that in
higher dimensions ETs are a dominating life form.

Theorems \ref{thm:dim_le_three} and \ref{thm:dim_four} follow from the
concrete computations of \cite{AlexeevErdal_dim4} where it is proved
that all periodic decompositions in dimension $\le3$ are Delaunay
(easy) and that in dimension $4$ there are exactly 2 non-Delaunay
decompositions (hard). Both of them are sub decompositions of the
unique maximal dicing with 9 hyperplanes, which I will denote by
$\Delta_{RT}$. The dimension of the stratum for $\Delta_{RT}$ in
$\Avor$ is~1. On the other hand, $h^0(\Delta,\bl)=2$ and so there must
be a second irreducible component.  The secondary polytope
$\Sigma(\Delta_{RT})$ is a square which corresponds to
$\bP^1\times\bP^1$, which gives $\bP^2$ after dividing by
$\Sym\Delta_{RT}=\bZ_2\times\bZ_2$.

The lower bound in theorem \ref{thm:dim_large} follows by computing a
particular example, the Delaunay decomposition of the classical
lattice $D_n$. One knows (see f.e. \cite{ConwaySloane93}) that there
are 3 maximal cells in this decomposition, 2 copies of a hemicube
$h\gamma_n$ (with $2^{g-1}$ vertices) and a crosspolytope $\beta_n$
(with $2g$ vertices), and that these polytopes are simplicial.
Applying formula~\ref{eqn:main} gives the bound.  The upper bound
follows from a simple observation that for a polytopal decomposition
the dimension of a secondary polytope $\Sigma(\Delta)$ is always less
than the volume of $|\Delta|$ (in the lattice units of volume). The
volume of $\XR/X$ is $g!$. If the decomposition is not polytopal but
is the pullback of a polytopal decomposition in $\bR^{g-a}$, one also
has to add the term $a(a+1)/2+a(g-a)$ to $(g-a)!$ to account for the
abelian part. However, it is easy to see that for $g\ge3$ the bound
$g!$ is greater.

I note another nice application of formula~\ref{eqn:main}: lattice
$E_8$ whose cells are 135 copies of a crosspolytope $\beta_8$ and 1920
copies of a simplex $\alpha_8$. The computation gives
$\dim\str[E_8]=h^0(E_8,\bl)=945$. This gives the largest ET in
dimension 8 that I am aware of. The normalization of the closure of
this stratum is the quotient of $(\bP^7)^{135}$ by a finite group.

\section{Dimension 4 case in detail}

The following is a classical description, due to Voronoi, of the 2nd
Voronoi decomposition in dimension 4. The picture below gives a
schematic view of a cross-section of this decomposition. First of all,
in dimension~4 $\tau_{\vor}$ is a subdecomposition of the perfect
decomposition. The perfect decomposition has cones of two types: 1st
domain, which are simplicial (with 10 sides) and cones with 64 sides.
In $\tau_{\vor}$ the cones of the second type are subdivided into 64
simplicial cones, 48 of which belong to the so called 2nd domain, and
16 -- to the 3rd domain. As we have indicated, the 2nd domain cones
have 1st domains as ``across-the-border'' neighbors, and 3rd domains
-- 3rd domains again.

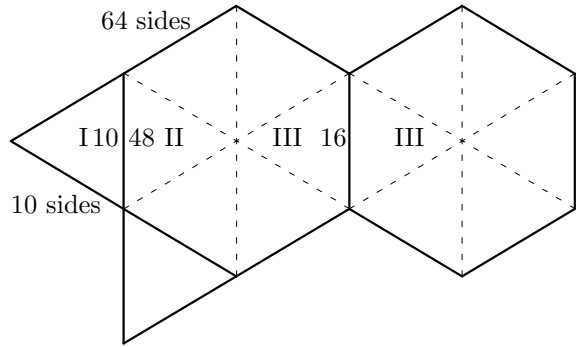
\begin{figure}[h]
  \begin{center}
    \setlength{\unitlength}{0.03cm}
    \begin{picture}(250,152)(0,-30)
      \thicklines
      \drawline(50,30)(0,60)(50,90)(50,30)
      \drawline(50,90)(100,120)(150,90)(150,30)(100,0)(50,30)
      \drawline(150,90)(200,120)(250,90)(250,30)(200,0)(150,30)
      \put(0,28){10 sides}
      \put(40,110){64 sides}
      \put(30,58){I}
      \put(68,58){II}
      \put(116,58){III}
      \put(170,58){III}
      \drawline(50,-30)(50,30)
      \drawline(50,-30)(100,0)
      \thinlines
      \dashline{3}(50,30)(150,90)
      \dashline{3}(50,90)(150,30)
      \dashline{3}(100,0)(100,120)
      \dashline{3}(150,30)(250,90)
      \dashline{3}(150,90)(250,30)
      \dashline{3}(200,0)(200,120)
      \put(52,58){48}
      \put(36,58){10}
      \put(137,58){16}
    \end{picture}
    \caption{2nd Voronoi decomposition in dimension 4}
  \end{center}
\end{figure}

It is exactly the face between two 3rd type domains where an ET
attaches itself. The way it happens is easier to see on a dual
picture. A dual picture is the picture for a polarization function on
$\tau_{\vor}$ ($\Avor$ is projective, \cite{Alexeev_CMAV}). Two
maximal cones of 3rd type on the dual picture correspond to points,
and the face between them -- to an interval. In the moduli space this
gives a $\bP^1$ (which we have to divide by automorphisms which gives
a $\bP^1$ again). If we denote by $\Delta_{RT}$ the corresponding
decomposition then $\dim\str[\Delta_{RT}]$ in $\Avor$ is one. One the other 
hand, explicitly this decomposition consists of 
\begin{enumerate}
\item a copy of a cyclic polytope $C_6$ with 6 vertices,
\item the inverse of it under the involution $x\mapsto -x$,
\item 18 simplexes.
\end{enumerate}
Hence, according to the formula of
example~\ref{exmp:decomps_with_simplicial_cells} the dimension of the
corresponding stratum in $\APg$ is $2\cdot l_{C_6}=2$. Hence, we must
have an ET. This ET on the dual picture is represented by a square
and is isomorphic to $\bP^1\times\bP^1$. After dividing by the
automorphism group $\bZ^2\times \bZ^2$ this gives a $\bP^2$ in the
coarse moduli space. As the picture suggests, the intersection of the
main component in $\APg$ and the ET is a diagonal in
$\bP^1\times\bP^1$.

\begin{figure}[h]
  \begin{center}
    \setlength{\unitlength}{0.05cm}
    \begin{picture}(80,60)(0,0)
      \thinlines
      \drawline(10,30)(40,60)(70,30)(40,0)(10,30)
      \thicklines
      \drawline(10,30)(70,30)
      \put(1,28){III}
      \put(71,28){III}
      \put(28,55){IV}
      \put(47,0){IV}
    \end{picture}
    \caption{Part of the dual picture with an ET attached}
  \end{center}
\end{figure}
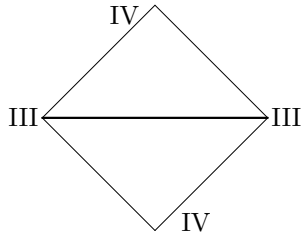

A vertex marked $IV$ corresponds to a non-Delaunay triangulation
obtained by triangulating one copy of $C_6$ in one of the two possible 
ways, and the other -- in the non-symmetric way. Clearly, this
triangulation cannot be Delaunay, indeed it does not even have the
basic symmetry $x\mapsto -x$. The same can be said about the
corresponding variety $P$.

\section{Relationship between ETs and the Jacobian locus}

The Torelli map $M_g\to A_g$ extends to a functorial morphism from the
Mumford-Deligne compactification $\overline{M}_g$ to $\APg$. On
$k$-points, it sends a stable curve $C$ to a triple $\Pic^0 C\acts
\Jac^{g-1}C \supset\Theta_{g-1}$, where $\Jac^{g-1}C$ is the 
moduli space of semistable rank 1 sheaves of degree $g-1$ on $C$.
\cite{Alexeev_CompJacobians} contains an algorithm for computing this
triple in terms of the dual graph of a stable curve $C$.

By checking all genus 4 stable curves we find that precisely 1 of them
maps to a point meeting an ET. It is a curve all of whose components
are $\bP^1$'s and the dual graph is the bipartite graph $K_{3,3}$. By
the Kuratowski theorem this is the minimal (in terms of genus) graph
that is not planar. The image of this curve in $\APg$ is precisely the 
center $1\in\bP^1$ connecting two 3 type domains on the above
figure. For any planar graph the corresponding Delaunay decomposition
belongs to the 1st domain (or one of its faces) and so is away from an 
ET. Therefore, we obtain

\begin{thm}\label{thm:planar}
  For a stable curve of genus $g\le4$ its Torelli image $[C]\in\APg$
  meets an ET if and only if the dual graph $C(\Gamma)$ is not
  planar. 
\end{thm}

One can easily make this theorem a little stronger: it is true for a
stable curve of arbitrary genus such that $h^1(C(\Gamma))\le4$. This
leads me to make the following
\begin{conj}
  Theorem \ref{thm:planar} holds in arbitrary genus.
\end{conj}

It feels that the non-planarity of $C$ somehow opens up new dimensions
for deformations for its jacobian.  Here is some additional
evidence in favor of this conjecture:
\begin{enumerate}
\item If $\Gamma(C)$ is planar, the corresponding Delaunay
  decomposition belongs to the simplest part of the 2nd Voronoi
  decomposition, the so called 1st domain.  If there is any part of
  $\Avor$ where ``nothing tricky happens'', that should be it.
\item It is plausible that if $\Gamma(C)$ contains a subgraph
  $\Gamma(C')$ for a curve $C'$ whose jacobian has extra deformations
  then the same must be true for the curve $C$ itself. By the
  Kuratowski theorem a graph $\Gamma$ is non-planar if and only if it
  contains either a $K_{3,3}$ or a $K_5$. The case of $K_5$,
  therefore, becomes a crucial test for the validity of the
  conjecture.
\end{enumerate}

\ifx\undefined\bysame
\newcommand{\bysame}{\leavevmode\hbox to3em{\hrulefill}\,}
\fi

\end{document}